\documentclass[12pt,a4paper]{article}
\usepackage[utf8]{inputenc}
\usepackage[english]{babel}

\usepackage{amsmath}
\usepackage{amssymb}
\usepackage{latexsym}
\usepackage{srcltx}
\usepackage{graphics}
\usepackage{color}
\usepackage{epsfig}

\newtheorem{theorem}{Theorem}

\newtheorem{cor}{Corollary}
\newtheorem{remark}{Remark}

\newcommand{\re}{{\mathbb R}}

\newcommand{\n}{{\mathbb N}}

\newcommand{\z}{{\mathbb Z}}

\newcommand{\bw}{{\mathbf{w}}}

\date{}

\author{M. Karapetyants 
\thanks{Moscow Institute of Physics and Technology, Regional Scientific and Educational Mathematical Center of Southern Federal University {e-mail: \tt\small
karapetyantsmk@gmail.com}},
V. Yu. Protasov
\thanks{University of L'Aquila, Moscow State University, Higher School of Economics {e-mail: \tt\small
v-protassov@yandex.ru}}   }

\title{The spaces of dyadic distributions
\thanks{
The second author is supported with Russian Foundation for Basic Research grants 
No 20-01-00469 and 19-04-01227.
}}

\begin{document}
\maketitle

\begin{abstract}

In this paper the spaces of distributions on a dyadic half-line, which is the positive half-line equipped with the digitwise binary addition and Lebesgue measure, are studied. We prove the non-existence of such a space of dyadic distributions which satisfies a number of natural requirements (for instance, the property of being invariant with respect to the Walsh-Fourier transform) and, in addition, is invariant with respect to multiplication by linear functions. This, in particular, allows the space of dyadic distributions suggested by S. Volosivets in 2009 to be optimal. We also show the applications of dyadic distributions to the theory of refinement equations as well as wavelets on a dyadic half-line.

\smallskip




\end{abstract}

\begin{center}
\textbf{1. Introduction}
\end{center}
\bigskip

The existence of more or less relevant space of dyadic distribution was discussed in the literature and in conferences since  early 2000s. The dyadic half-line is the positive half-line equipped with the specific addition operation, namely, if $x \, = \, \sum_{k\in \z} x_k2^{k}, \, y \, = \, \sum_{k\in \z} y_k2^k$ are arbitrary positive numbers presented in their binary form (each series starts with the sequence of zeros), then their dyadic sum is defined as follows:
$x \oplus y\, = \, \sum_{k\in \z} 2^{k}\,(x_k+y_k)_2$, where $(x_k+y_k)_2$ is the sum modulo 2 of $k$-th elements in the binary expansion of $x$ and $y$ respectively. Herewith, $x \ominus y = x \oplus y$. The Lebesgue measure on the dyadic half-line $\re_+$ coincides with its standard analogue, and so do all the spaces~$L_p(\re_+)$. The Walsh functions play the role of exponents in the dyadic harmonic analysis. The Walsh-Fourier transform is isometric in the space~$L_2(\re_+)$. More about the properties of functions on a dyadic half-line and its applications can be found in~\cite{GES, SWS}.

In 2007 B. Golubov \cite{G} defined a space of distributions on a dyadic half-line~$\re_+$. It consists of continuous linear functionals on the space of test functions $ D_d(\mathbb R_+)$, which is the space of infinitely differentiable (in the dyadic sense) functions on $\mathbb{R}_+$ such that \\

1)the supports of each function $ \varphi $ alongside with all its dyadic derivatives ${\varphi}^{\alpha}, \ \alpha \in \mathbb N, $ are contained in some dyadic interval $ \delta $; \\

2)for all $ \alpha \in \mathbb N $, the sequence $ \psi \ast ( { {\Lambda}^{\alpha} }_n )(x) $ uniformly converges to $ {\psi}^{\alpha} $ on $ \mathbb R_+ $ as $ n \to \infty $, where

$$
{ {\Lambda}^{\alpha} }_n = \int_{0}^{2^n} { (h(t)) }^{- \alpha} \psi(x, t) dt, \ x \in \mathbb R,
$$

$$
h(x) = 2^{-n}, \ 2^n \leq x < 2^{n+1}, \ n \in \mathbb Z,
$$

$ \psi(x, t) $ is the the generalized Walsh functions (see the corresponding definition in section 2). \\

The space of rapidly decreasing in the neighbourhood of infinity functions $ S_d(\mathbb R_+) $ was also suggested. It consists of infinitely smooth functions (in dyadic sense) such that for each $ \alpha, \ \beta \ \in { \mathbb Z }_+$, we have 

$$
\lim_{x \to \infty} {(h(x))}^{- \beta} { \varphi }^{ \alpha }(x) \ = \ 0.
$$

In this case the space of distributions was also defined as the space of linear continuous functionals on $ S_d(\mathbb R_+) $.

The main problem of the spaces $ D_d$ and $ S_d$ is that they are not invariant with 
 respect to the Walsh-Fourier transform. Thus, the Walsh-Fourier transform cannot  be well defined on the corresponding spaces of distributions~$D_d', S_d'$. This problem, however, was solved by S.S. Volosivets in 2009. Note that similar constructions in the space $L_2$ on other groups appeared in~\cite{VVZ, B, T}. In~\cite{V} he suggested another, more narrow than $D_d$, space of test functions $H_d$, which is the space of ''dyadic-analytic'' functions $f$ which are the finite linear combinations of the indicator function $\chi_{\Delta_{j,k}}$, where $\Delta_{j,k} \, = \, [2^{-j}k, 2^{-j}(k+1))$ is the dyadic interval of the rank~$j\in \mathbb{Z}$. For all such functions $f$ it holds that $f(\cdot + h) - f(\cdot) \equiv 0$ for an arbitrary $h \in (0, 2^{-n})$, where $n$ is the highest rank of the intervals in the linear combination.

The topology is defined by convergence to zero: $f_k \to 0$ as $ k \to \infty $ if the ranks of the intervals of the sequence $\{f_k\}_{k \in \n}$ are bounded above and are contained in the fixed segment, and the sequence itself converges to zero pointwise. It is easy to show that $H_d$ is a complete linear space invariant with the respect to the Walsh-Fourier transform. Linearity and completeness can be checked directly and the last property follows from the fact that if ${\rm supp} f\,  \subset \, [0, 2^m]$ and $n$ is the highest rank of the intervals in~$f$, then ${\rm supp} \widehat f\,  \subset \, [0, 2^n]$ and $m$ is the highest rank of the intervals in~$\widehat f$. Thus the Walsh-Fourier transform maps the space into itself.

 The Walsh-Fourier transform on the corresponding space of distributions $H_d$ is defined as usual: for each $f\in H_d'$, we have $(\widehat f, \varphi) \, = \, (f, \check \varphi), \, \varphi \in H_d$, where $\check \varphi$ means the inverse Walsh-Fourier transform.

\bigskip

\begin{center}
\textbf{2. On the possibility of multiplication by smooth functions. \\ The main result.}
\end{center}
\bigskip

The space of distributions suggested by S. Volosivets is very convenient due to its simplicity and variety of applications. Since the space of the test functions is very narrow (it only contains functions generated by binary dilates  and shifts of the function $\chi_{[0,1)}$), the space of distributions is rather wide. For instance, every locally summable on $\re_+$ function~$f$ belongs to $H_d'$ with its Walsh-Fourier transform. The exponent also belongs to $H_d'$ which is not true for the classical Schwartz space.

The main disadvantage of $H_d'$ is that it is not invariant with  respect to multiplication by smooth functions, in particular, by polynomials. For example, not for each function $f \in H_d'$ is true that also $xf \in H_d'$. The natural question is whether there is a space of distributions on the dyadic half-line that is invariant with respect to both the Walsh-Fourier transform and the multiplication by smooth functions, e.g. by polynomials? Theorem \ref{th.10} gives a negative answer and establishes the non-improvability of the space~$H_d'$. There is no extension of $H_d'$ that allows the multiplication even by linear functions.


First we need to introduce some notation. For arbitrary $x, y \in \re_+$, we denote 
$(y, x) \, = \, \sum_{k \in \z}^{\infty} y_k x_{-1-k}$, where $x_i, y_i$ are the digits in the binary expansion of $x$ and $y$ respectively. This sum always contains a finite number of nonzero terms.
For an integer $k \ge 0$, the Walsh function is defined as $\bw_{k}(x) \, = \, (-1)^{(k, x)}$ and $ \psi(x, y) \ = \ \bw_{[y]}(x) \cdot \bw_{[x]}(y) $, where $ [y] $ is the integer part of $ y $. The Walsh-Fourier transform of the function $f\in L_1(\re_+)$ is $\widehat f (y)\, = \, \int_{\re_+} \psi(x, y) f(x) \, dx $, and it can be extended to $L_2(\re_+)$ in a usual way. The Walsh-Fourier transform is an invertible orthogonal transform of $L_2(\re_+)$ \cite{GES, SWS}.

What can be the space of the test functions  to define distributions on a dyadic half-line?
 It is natural to require  that this space  contains the indicator function $\chi_{[0,1)}$
  and is invariant with respect to integer shifts of a function $f(x) \, \mapsto f(x\oplus 1)$ as well as with respect to binary contraction and expansion $f(x) \, \mapsto f(2x), \ f(x) \, \mapsto f(x/2)$. In this case it already contains $H_d$. Thus, $H_d$ is the smallest 
  by inclusion  functional space satisfying these requirements. The property of being invariant with respect to the Walsh-Fourier transform is fulfilled automatically. Indeed, $\widehat \chi_{[0,1)}\, = \, \chi_{[0,1)}$, so the Walsh-Fourier transform maps $H_d$ into itself. The question arises whether it is possible to extend $H_d$ so that it would be also invariant with  respect to the multiplication by algebraic polynomials? If so, the corresponding space of distributions would also possess this property: if $f\in H_d'$, then $xf$ is defined as in the classical case: $(xf, \varphi) \, = \, (f, x\varphi), \, \varphi \in H_d$.

Theorem \ref{th.10} provides a negative answer to the question above under another natural condition: the space of distributions must contain the functional space $L_2(\re_+)$, which means that each $g \in L_2(\re_+)$ acts naturally on
$H_d$, e.g., the integral~$\int_{\re_+} fg \, dx$ is defined for each function $f \in H_d$.

\begin{theorem}\label{th.10}
There is no space of measurable functions on a dyadic half-line that contains the indicator function~$\chi_{[0,1)}$, is invariant with respect to both the Walsh-Fourier transform and the multiplication by linear functions, and on which every element of the space $ L_2(\re_+) $  
acts by the formula of the inner product.
\end{theorem}

{\tt Proof}. Suppose such a space exists, we denote it by $\tilde H_d$. Computing the Walsh-Fourier transform of the function $f(x) = x\chi_{[0,1)}(x)$, we obtain:
$$
\widehat {f}(y)\ = \
\int_0^1 x \cdot \psi(x, y) dx \ = \ \int_0^1 x \cdot \bw_{[x]}(y)
\cdot \bw_{[y]}(x) dx \ = \ \int_0^1 x \cdot \bw_{[y]}(x) dx\, ,
$$
since $ \bw_{[x]}(y) \ = \ \bw_{0}(y)\, = \, 1 $.
This integral depends only on the integer part of $ y $.
We calculate it for $y \in [2^n, 2^{n}+1)$, when $[y] = 2^n, \,
n \in \n$. So, $[y] = 10\ldots 0$ ($n$ zeroes);
then $\bw_{[y]}(x) \, = \, (-1)^{(2^n, x)} \, = \, (-1)^{x_{-n-1}}$.
Thus, if $y\in [2^n, 2^n+1)$, we obtain
$$
\int_0^1 x \cdot \bw_{[y]}(x) dx \ = \ \int_0^1 x \cdot x_{-n-1} \,  dx \ = \
\sum_{k=0}^{2^{n+1}}(-1)^k\int_{k2^{-n-1}}^{(k+ 1)2^{-n-1}}\, x  \,  dx\ =  \
$$
$$
\ \sum_{k=0}^{2^{n+1}}(-1)^k\frac{(k+ 1)^2-k^2}{2^{2n+3}} \ = \ - 2^{-(n+1)}\, .
$$
Finally,
$$
\widehat {f}(y)\quad  =  \quad  - \, 2^{-(n+1)} \ , \qquad
y \in {[2^n, 2^{n}+1)}\,  , \quad n \in \n \, .
$$
By the assumption,  $\widehat{f}(y)  \in \tilde H_d$, therefore
$y\widehat{f}(y)  \in \tilde H_d$. On the other hand, each element from $L_2(\re_+)$ acts naturally on $\tilde H_d$. We choose the following function $g\in L_2(\re_+)$:
$$
    g(y)  \ = \ \left\{ \begin{array}{cl}
    -\, \frac{1}{n+1}, & \ y \in [2^n; 2^n + 1), \ n \in \mathbb{N} \, , \\
    0, & \ \mbox{ else}.
    \end{array}\right.
$$
Then
$$
\Bigl(g\, , \,  y \, \widehat{f}(y)\Bigr) \ =\
- \, \sum_{n = 1}^{\infty}\, \frac{1}{n+1} \, \int_{2^n}^{2^n+1} \, \frac{y}{2^{n+1}}\, dy\ = \ - \, \sum_{n = 1}^{\infty}\, \frac{1}{n+1} \, \frac{y^2}{2^{n+2}}\Bigl|_{2^n}^{2^n+1} \ = \
$$
$$
\ = \ - \, \sum_{n = 1}^{\infty}\, \frac{1}{n+1}\, \frac{2^{2n} \, - \, (2^n+1)^2}{2^{n+2}}\ =\
\sum_{n = 1}^{\infty}\,  \frac{\frac12 + 2^{-n-2}}{n+1}\ = \ \infty.
$$
Consequently, $\Bigl(g\, ,  \, y \widehat{f}(y)\Bigr)$  is not defined, which leads to a contradiction.

{\hfill $\Box$}
\smallskip

\begin{remark}\label{r.20}{\em
One may also define a space of dyadic smooth rapidly decreasing functions $Q_d$. Dyadic smoothness means that for a function $f$, the following holds:
$$
     \Bigl\| f(x \oplus t)\,  - \, f(x) \Bigr\|_{2} \  \leq \ C(\alpha) \cdot {  t  }^\alpha\, \qquad
    t > 0,\
$$
for each $ \alpha >0$, and $f$ is rapidly decreasing if $ |f(x)| \leq C{(x+1)}^{-n} $, for each $ n \in \mathbb N $, $ x \in \mathbb R_+ $.
The space $ Q_d $ is invariant with respect to the Walsh-Fourier transform, but it is also invariant with  respect to multiplication by dyadic smooth rapidly decreasing functions. However, the dyadic smoothness is not the same as the smoothness in the classical sense, that is why $ Q_d $ is not invariant with respect to multiplication by linear functions.
}
\end{remark}

\bigskip

\begin{center}
\textbf{3. Applications to the wavelet theory}
\end{center}
\bigskip

The first examples of systems of wavelets on the dyadic half-line and its various generalizations can be found in the paper of Lang~\cite{L}, more general constructions were presented in~\cite{F, PF, RF}. The wavelets on the Abelian groups were also studied~\cite{LBK}.
To obtain a system of  wavelets one needs to solve the {\em refinement equation}, which is a functional equation on a function $ \varphi $ with binary expansion of the argument $x$ :
$$
\varphi \ = \ \sum_{k=0}^{2^n} c_k \varphi(2x \ominus k), \ x \in \re_+
$$
Such equations are also used when studying dyadic approximate algorithms~\cite{K}. The theory of refinement equations on the classical real line~$\re$ was developed  by the end of 1980s~\cite{D,NPS}. However, many questions remained unsolved on a dyadic half-line. Does the refinement equation always have a solution and, if so, which class does the solution belong to? Will the solution be unique up to multiplication by a constant? If the solution is
compactly supported,  what is the length of its support? The space of dyadic distributions~$H_d$ allows us to fully answer these questions in Theorem~\ref{th.30}. First let us introduce some further notation.

The solution of the refinement equation is called a refinable function, which is a fixed point 
of the {\em transition operator $T$}:
$Tf \ = \ \sum_{k=0}^{2^n} c_k f(2x \ominus k)$.
Set $m(y) = \frac12 \, \sum_{k=0}^{2^n} c_k \bw_k (y)$. This Walsh polynomial is called the {\em mask} of the refinement equation. It is known~\cite{PF} that the studying of the general refinement equations could be reduced to the case $\sum_k c_k = 2$, which is equivalent to $m(0)=1$.

\begin{theorem}\label{th.30}
For each sequence of complex coefficients $\{c_k\}_{k=0}^{2^n}$, the sum of which equals $2$, the refinement equation has a unique up to multiplication by a constant solution in terms of distributions $\varphi \in H_d'$. The support of this function $\varphi$ lies in the segment $[0, 2^n]$ and the Walsh-Fourier transform of $\varphi$ is given by the formula:
$$
\widehat \varphi (y) \ = \ \prod_{j=1}^{\infty} m\bigl(2^{-j}y \bigr)
$$
Moreover, for each finite summable function $f\in  L_2(\re_+)$,
the sequence $T^kf$ converges in $H_d'$ to the solution of the refinement equation $c\, \varphi$, where $c = \int_{\re_+}f(x)dx$.
\end{theorem}
{\tt Proof.} Using the properties of the Walsh-Fourier transform, we obtain: $\widehat {Tf}(y) \, = \, m\bigl(\frac{y}{2} \bigr)\widehat {f}\bigl(\frac{y}{2} \bigr)$. Consequently, for each $k$ we have
\begin{equation}\label{eq.prod}
\widehat {T^kf} (y) \ = \ \widehat {f}\bigl(2^{-k} y\bigr)\, \prod_{j=1}^{\infty}
m\bigl(2^{-j}y \bigr)\,  .
\end{equation}

We show that for each finite function $f\in L_1(\re_+)$, 
the product~(\ref{eq.prod}) converges uniformly on each segment $[0, 2^N]$.
Note that the function $m(y)$ is a Walsh polynomial of the power $2^n$, so it is constant on the dyadic intervals of rank~$n$. Since $\sum_{i}c_i = 2$, we deduce
$m(0) = 1$. Thus, $m(z) = 1$ for each $z \in [0, 2^{-n})$.
Therefore, as $k > n + N$, we have $2^{-k} y \in [0, 2^{-n})$ for each $y \in [0, 2^N]$, and hence $m\bigl(2^{-k} y\bigr) \, = \, 1$. So, on a segment~$[0, 2^N]$ each term in~(\ref{eq.prod}) with the number~$k >  n + N$ is identically equal to one on~$[0, 2^{N}]$, and the product converges on this segment. Thus, the product~(\ref{eq.prod}) converges uniformly on each compact set in~$\re_+$, consequently, it converges in the space of distributions in~$H_d'$. Therefore $T^kf$ converges in~$H_d'$ to a distribution~$\psi$. Then $T\psi = \psi$, e.g.,  $\psi$ is the solution of the refinement equation.

Since $f$ is finite, let its support be in $[0, 2^{\ell})$.
Then on the segment $[0, 2^{-\ell})$ the function $\widehat f$ is identically equal to $\widehat f(0) = \int_{\re_+} f\, dt\, = \, c$.
As $k > \ell + N$ we obtain $2^{-k} y \in [0, 2^{-\ell})$ for each $y \in [0, 2^N]$, hence $\widehat {f}\bigl(2^{-k} y\bigr) \, = \, c$. Thus, the product~(\ref{eq.prod}) converges as $k \to \infty$ to~$\, c\, \prod_{j=1}^{\infty} m\bigl(2^{-j}y \bigr)$.
We see that  the inverse Walsh-Fourier transform of this product is nothing but a solution of the refinement equation. Therefore, the inverse Walsh-Fourier transform of~$\prod_{j=1}^{\infty} m\bigl(2^{-j}y \bigr)$ is also a refinable function, which is denoted by~$\varphi$. Finally, for each finite function $f$, the sequence $T^kf$ converges to $ \varphi $.

Now let $\varphi_0$ be an arbitrary finite solution of the refinement equation. Since $T\varphi_0 = \varphi_0$, it follows that $T^k\varphi_0 = \varphi_0$ for each~$k$. Hence, the sequence $T^k\varphi_0$ converges to $\varphi_0$, the function $\varphi_0$ is proportional to $\varphi$ (namely, $\varphi_0 = \varphi \int_{\re}\varphi_0 dx$). From this the uniqueness of the solution follows.
Finally, to prove that the support of the solution lies in
$[0, 2^n]$, it is enough to consider an arbitrary function $f$ supported by this segment and apply the operator $T$.
The function $Tf$ is also supported on the same segment and so do all the functions $T^kf$, therefore the limit (which is the solution of a refinement equation) of the sequence $T^kf$ is also supported on $[0, 2^n]$.

{\hfill $\Box$}
\smallskip

Theorem~\ref{th.30} provides tools for studying the properties of many refinements equations. For instance, the following fact is useful in probability theory, approximation theory and the theory of subdivision schemes.
\begin{cor}\label{c.10}
If all the coefficients $c_k$ of a refinement equation are non-negative, then its solution  $\varphi$ normalized by the condition $\int_{\re_+}\varphi dx = 1$, is a non-negative distribution.
\end{cor}
{\tt Proof.} Take $f=\chi_{[0,1)}$ and consider the sequence $T^kf$. Since the operator $T$ respects the non-negativity of functions, it follows that all the elements of that sequence are also non-negative. Hence, its limit, which is the solution, is non-negative as well.

\bigskip
\noindent \textbf{Acknowledgments.} The authors are grateful to S.S. Volosivets for valuable remarks and interesting  discussions and to an anonymous referee for a thorough reading and 
many useful advises.
The second author was supported by the Program of fundamental research of Higher School of Economics (National Research University) and was financed  within the framework of the state support for the leading universities of the Russian Federation ''5-100''.

{\hfill $\Box$}

\end{document}